\documentclass[twoside,11pt]{article}

\usepackage[utf8]{inputenc}
\usepackage[T1]{fontenc}
\usepackage{amsmath,amssymb,amsthm}
\usepackage{graphicx}
\usepackage{url}

\theoremstyle{plain}
\newtheorem{thm}{Theorem}[section]

\newtheorem{lem}[thm]{Lemma}

\newtheorem{defn}[thm]{Definition}

\newtheorem{ex}[thm]{Example}
\newtheorem{hyp}[thm]{Hypothesis}
\newtheorem{stmt}[thm]{Statement}

\begin{document}

\title{On the Connected Components of Structure Graphs of Associative Rings}

\author{%
Andrey M. Chernev\thanks{Corresponding author. Email: chernevAM@mpei.ru}\\
National Research University ``Moscow Power Engineering Institute'', Moscow, Russia
\and
Nikolay A. Chernev\thanks{Email: nachernev@gmail.com}\\
Bauman Moscow State Technical University, Moscow, Russia
}

\date{}
\maketitle

\begin{abstract}
The article studies the notion of the structure graph of an associative ring with unity.
It is known that the structural graph of a ring with $n$ elements has at least $n$ connected components.
The article investigates when this lower bound is attained. A hypothesis is stated that the only rings
of order $n$ having exactly $n$ connected components are the rings $\mathbb{Z}_2$, $\mathbb{Z}_3$, $\mathbb{Z}_7$, $GF(4)$, and their direct sums.
The hypothesis is proved for all rings except residue rings modulo a prime of the form $7+12k$, $k\geq181$,
and Galois fields of characteristic~2 with extension degree $\geq17$.
\end{abstract}



\markboth{Andrey M. Chernev, Nikolay A. Chernev}{On the Connected Components of Structure Graphs of Associative Rings}

\section{Introduction}

The subject of study in this article is the directed graph defined on the Cartesian square of an associative ring, in which an edge connects a pair of elements of the ring with the results of the ring operations -- addition and multiplication. This construction was proposed by A.~Lipkovski in \cite{link_Lipkovski_1995}. This topic was further developed in \cite{link_Lipkovski_2012}, \cite{link_AS_2013}, \cite{link_HS_2013}, \cite{link_Lipkovski_2015},
\cite{link_Lipkovski_2017} (the article \cite{link_Lipkovski_2017} is English version of the \cite{link_Lipkovski_2015}). Among recent works, we can mention \cite{link_DSL_2022} and \cite{link_Matovich_2025}.

In the English-language literature, this object is usually called a \emph{digraph} (or \emph{directed graph}) of ring.
We find this term too general and prefer \emph{structure graph}, the name used by A.~Lipkovski in his Russian-language papers.

\section{Preliminary Information}

Let us recall some definitions and facts used in this article that relate to graph theory, number theory, and algebra.

A directed graph is called \emph{functional} if the out-degree of every vertex is $1$.

If on some set $X$ a function $f: x \to f(x)$ is given, then one can associate to this function a directed graph with the set of vertices is $X$ and the set of edges is the pairs $(x, f(x))$. We will call such a graph the \emph{graph of the mapping $f$}. Of course, it will be functional.

We will call a directed graph \emph{connected (weakly connected)} if its underlying undirected graph (with the same sets of vertices and edges) is connected.

A \emph{(weakly) connected component} of a directed graph $G$ is a maximal (weakly) connected subgraph of $G$.

Any weakly connected component of a functional graph consists of a directed cycle, some vertices of which are roots of inverted trees (arrows go from leaves to roots).

The \emph{Legendre symbol} $\left(\frac{a}{p}\right)$ for a prime $p$ and $a \not\equiv 0 \pmod{p}$ is defined as follows:
\[
\left(\frac{a}{p}\right) =
\begin{cases}
    1, & \text{if } x^2 \equiv a \pmod{p} \text{ has a solution,} \\
    -1, & \text{if } x^2 \equiv a \pmod{p} \text{ has no solution.}
\end{cases}
\]

The following assertions are valid (see, for example, \cite{link_Vino}):
\[
\left(\frac{ab}{p}\right) = \left(\frac{a}{p}\right)\left(\frac{b}{p}\right),
\]
\[
\left(\frac{-1}{p}\right) = (-1)^{\frac{p-1}{2}},
\]
\[
\left(\frac{q}{p}\right) = (-1)^{\frac{p-1}{2}\frac{q-1}{2}} \left(\frac{p}{q}\right).
\]

A \emph{Galois field} ($GF(n)$) is any field of finite order. The order of a Galois field is always equal to $p^n$ for some prime $p$.

Galois fields of the same order are isomorphic and are obtained from $\mathbb{Z}_p$ by adjoining a root of an irreducible polynomial of degree $n$ over
$\mathbb{Z}_p$. Furthermore, $GF(p^k) \subseteq GF(p^m)$ if and only if $k \mid m$.

Let us dwell in more detail on the definitions and assertions of the structure theory of rings.

\begin{defn}
The \emph{Jacobson radical} of an associative ring is the intersection of all its maximal right ideals.
\end{defn}

The intersection of maximal left ideals also coincides with the Jacobson radical.

\begin{defn}
An associative ring is called \emph{semisimple} if its Jacobson radical is zero.
\end{defn}

\begin{defn}
An associative ring is called \emph{right Artinian} if any nonempty set of its right ideals has a minimal element.
\end{defn}

In particular, any finite ring is right Artinian.

\begin{stmt}\label{stmt:nilpotent}
The Jacobson radical of a right Artinian ring is nilpotent.
\end{stmt}
\begin{proof}
See Theorem~1.3.1 in \cite{link_Herst}.
\end{proof}

Let us present the main fact of the classification of semisimple finite rings.

\begin{thm}\label{thm:semisimple}
Any finite semisimple ring is isomorphic to a direct sum of matrix rings over some fields.
\end{thm}
\begin{proof}
This assertion is a consequence of several theorems from \cite{link_Herst}, although it is not contained there in explicit form. By the Wedderburn--Artin theorem (2.1.7 from \cite{link_Herst}), a semisimple right Artinian ring is isomorphic to a direct sum of matrix rings over division rings. Any finite ring is right Artinian. Any finite division ring is a field (Theorem~3.1.1, \cite{link_Herst}).
\end{proof}

\section{Structure Graphs}

\begin{defn}
The \emph{structure graph} of an associative ring $R$ is a directed graph whose vertices are pairs of elements of $R$, and whose edges
are directed from a vertex $(a, b)$ to the vertex $(a+b, ab)$.
\end{defn}

Previously, the notion of a structure graph was defined only for commutative rings, but this notion can be extended to noncommutative rings.

Let us present several assertions whose proofs follow directly from the definition of a structure graph.

\begin{stmt}\label{stmt:isomorphic}
If rings are isomorphic, then their structure graphs are isomorphic.
\end{stmt}

The converse is false. In \cite{link_HS_2013}, an example of non-isomorphic rings with isomorphic structure graphs ($\mathbb{Z}_4$ and $\mathbb{Z}_2[i]$) is given.

\begin{stmt}\label{stmt:subring}
If a ring $R$ is a subring of $S$, then the structure graph of $R$ is a subgraph of the structure graph of $S$.
\end{stmt}

\begin{stmt}\label{stmt:directsum}
If a ring $R$ is a direct sum of rings $A$ and $B$, then the structure graph of $R$ is the Cartesian product of the structure graphs of $A$ and $B$.
\end{stmt}

One of the questions studied in \cite{link_Lipkovski_2012} and \cite{link_Lipkovski_2015} was the study of the number of connected components of structure graphs.

Since the structure graph of a ring is a functional graph, each of its connected components ends with a cycle. Cycles of length~1 (loops) have a simple description.

\begin{stmt}\label{stmt:ncomponents}
A ring $R$ with $n$ elements contains exactly $n$ connected components ending with a loop (a cycle of length~1) $(a,0) \to (a,0)$ for any element $a \in R$ of the ring.
\end{stmt}

\section{The Main Hypothesis}

Thus, the structure graph of a ring with $n$ elements has at least $n$ connected components ending with a loop $(a,0) \to (a,0)$. A natural question arises: in what cases will there be no other connected components in the structure graph (an alternative formulation: in what cases does the structure graph of a ring contain no cycles of length $>1$)?

In \cite{link_Lipkovski_2012}, conditions following from the existence of cycles of length 2, 3, and 4 in the structure graph were studied. In particular, a criterion for the existence of cycles of length~2 was given.

\begin{stmt}\label{stmt:cycle2}
The structure graph of a ring $R$ contains a cycle of length~2 if and only if $R$ has nontrivial nilpotent elements.
\end{stmt}
\begin{proof}
See \cite{link_Lipkovski_2012}.
\end{proof}

However, the question in a more general formulation was not posed previously. The present work gives an almost complete answer to it (with the exception of some cases).

\begin{hyp}\label{hyp:main}
The structure graph of a ring whose order is equal to $n$ has exactly $n$ connected components if and only if the ring is a direct sum of rings, each summand of which is isomorphic to one of the rings $\mathbb{Z}_2$, $\mathbb{Z}_3$, $\mathbb{Z}_7$ (residue rings modulo the corresponding modulus), $GF(4)$ (the Galois field of 4 elements).
\end{hyp}

\begin{figure}[htbp]
    \centering
    \includegraphics[width=0.8\textwidth]{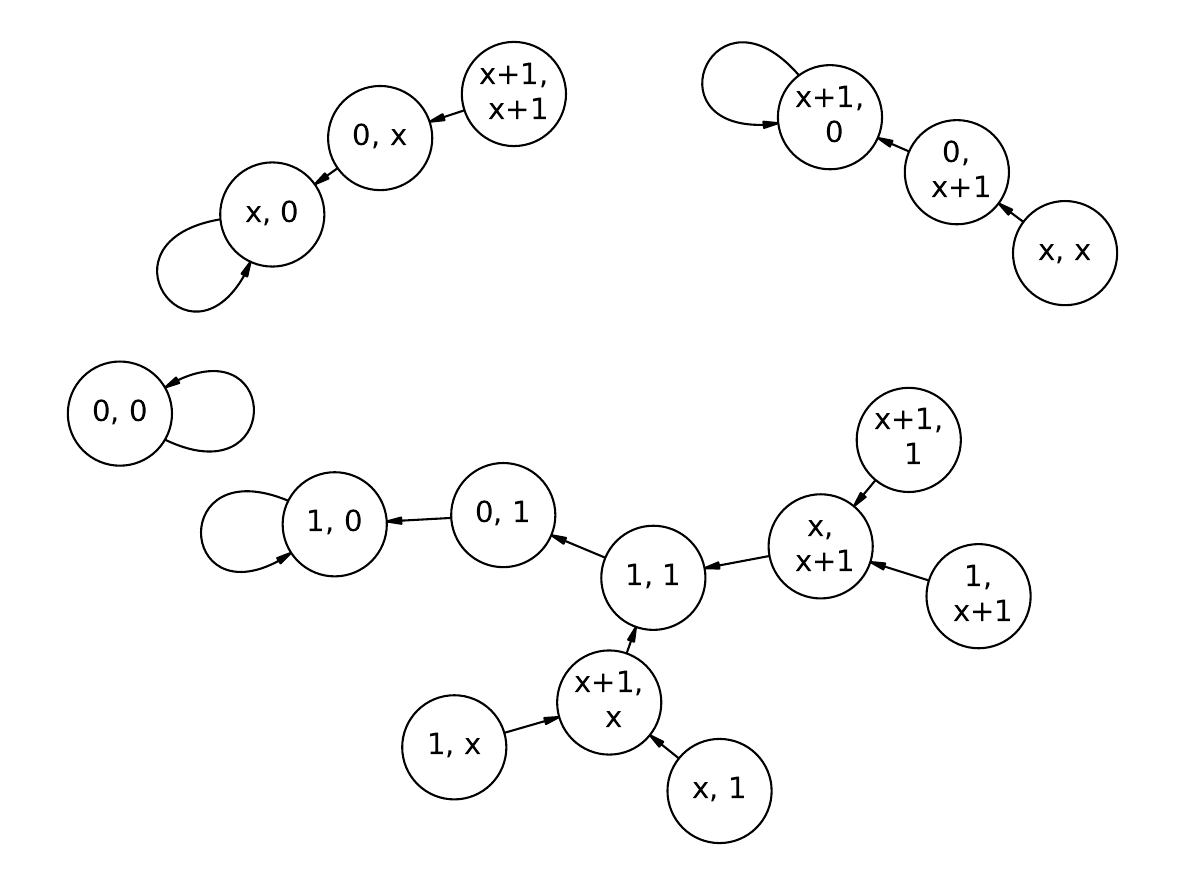}
    \caption{Structure graph of the ring $GF(4)$}
    \label{fig:GF4}
\end{figure}

In one direction, this assertion is proved by direct construction of the structure graphs of the corresponding rings
(for verification of computations see Section \ref{sec:verif} of this article), and by reference to Statement~\ref{stmt:directsum}. For example, Figure~\ref{fig:GF4} shows the structure graph of $GF(4)$ (the symbol $x$ denotes a root of the polynomial $x^2+x+1$).

The substantive assertion proved in this article is that in the structure graphs of any other rings there is a cycle of length $>1$.

In accordance with the classification theorem for finite semisimple rings (Theorem~\ref{thm:semisimple}), the proof of the existence of a cycle of length $>1$ in the structure graph of any ring other than a direct sum of $\mathbb{Z}_2$, $\mathbb{Z}_3$, $\mathbb{Z}_7$, $GF(4)$ splits into the following stages:

\begin{enumerate}
    \item Proof for $\mathbb{Z}_p$ ($p$ is a prime number, $p \neq 2, 3, 7$).
    \item Proof for Galois fields of extension degree $>1$ other than $GF(4)$.
    \item Proof for matrix rings of size $>1$ over finite fields.
    \item Proof for direct sums in which at least one of the summands is isomorphic to rings from (1)--(3).
    \item Proof for rings with a nonzero Jacobson radical.
\end{enumerate}

\section{Cycles of Length $>1$ for the Structure Graph of $\mathbb{Z}_p$ ($p$ Prime, $p \neq 2, 3, 7$)}

For rings of characteristic $\neq 2$, the problem can be reduced to the consideration of the graph of a mapping defined on the elements of the ring.

\begin{lem}\label{lem:g}
If in the graph of the mapping $g(x) = x - x^2$ defined on a ring $R$ there is a cycle of length $m > 1$, then in the structure graph of $R$ there is a cycle of length $2m$.
\end{lem}
\begin{proof}
Consider in the structure graph of the ring the path starting at the point $(-1, a)$:
\[
(-1, a) \to (a-1, -a) \to (-1, a - a^2) = (-1, g(a)).
\]
The first component repeats after 2 iterations, and the second one maps $a \to g(a)$. If $a = g^m(a)$ for some $m > 1$, then after $2m$ iterations the element $(-1, a)$ goes into itself, i.e., there is a cycle of length $2m$.
\end{proof}

Thus, the problem of proving the existence of cycles of length $>1$ in the structure graph of a ring reduces to proving the existence of cycles in the graph of the mapping $g(x) = x - x^2$ (in what follows, by $g(x)$ we mean this mapping).

\begin{figure}[htbp]
    \centering
    \includegraphics[width=0.6\textwidth]{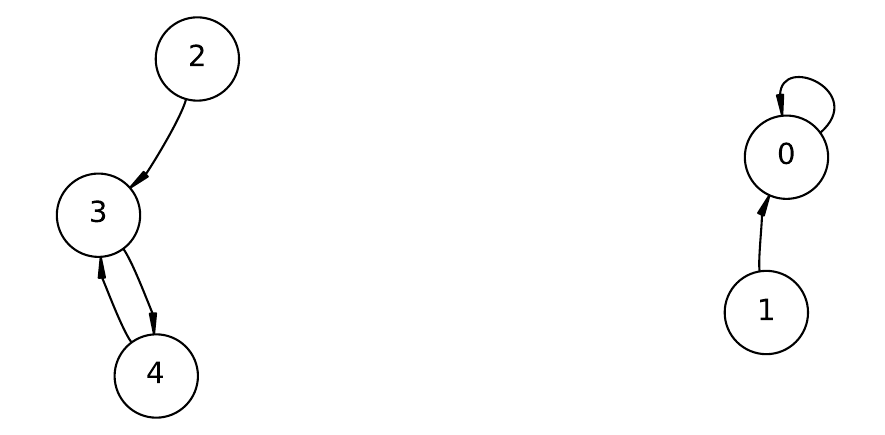}
    \caption{Graph of the mapping $g(x) = x - x^2$ in the ring $\mathbb{Z}_5$}
    \label{fig:Z5}
\end{figure}

\begin{ex}\label{ex:Z5}
As shown in Figure~\ref{fig:Z5}, in the ring $\mathbb{Z}_5$ the graph of the mapping $g(x)$ contains the cycle $(3) \to (4) \to (3)$. Therefore, by Lemma~\ref{lem:g}, the structure graph of $\mathbb{Z}_5$ contains a cycle of length~4:
\[
(4, 3) \to (2, 2) \to (4, 4) \to (3, 1) \to (4, 3).
\]
\end{ex}

Note that the graph of the mapping $g(x)$ always has only one loop $0 \to 0$. Therefore, the presence of at least 2 connected components for this graph entails the existence of a cycle of length $>1$.

\begin{thm}\label{thm:p1mod4}
If $p \equiv 1 \pmod{4}$, then the structure graph of the ring $\mathbb{Z}_p$ contains a cycle of length~4.
\end{thm}
\begin{proof}
Consider the equation $g(g(x)) = x$. For $x \neq 0$ it is equivalent to $(x-1)^2 = -1$. We have
\[
\left(\frac{-1}{p}\right) = (-1)^{\frac{p-1}{2}} = 1.
\]
Consequently, this equation has a solution. Thus, by Lemma~\ref{lem:g}, the structure graph contains a cycle of length~4.
\end{proof}

\begin{thm}\label{thm:p11mod12}
If $p \equiv 11 \pmod{12}$, then the structure graph of the ring $\mathbb{Z}_p$ contains a cycle of length $>1$.
\end{thm}
\begin{proof}
Consider the connected component of the graph of the mapping $g(x)$ containing $0$. The preimage of $0$ is the set $\{0, 1\}$. The preimage of $1$ must be a solution of the equation $x^2 - x + 1 = 0$. But for $p \equiv 11 \pmod{12}$ this equation has no solutions, since its discriminant is $-3$, and $-3$ is a quadratic non-residue for $p \equiv 11 \pmod{12}$. Indeed,
\[
\left(\frac{-3}{p}\right) = \left(\frac{-1}{p}\right)\left(\frac{3}{p}\right) = (-1)^{\frac{p-1}{2}}\left(\frac{3}{p}\right) = -\left(\frac{3}{p}\right),
\]
\[
\left(\frac{3}{p}\right)\left(\frac{p}{3}\right) = (-1)^{\frac{3-1}{2}\frac{p-1}{2}} = -1.
\]
Therefore
\[
\left(\frac{-3}{p}\right) = \left(\frac{p}{3}\right) = \left(\frac{2}{3}\right) = -1.
\]
Thus, no element other than $1$ and $0$ can be mapped to $0$, and, consequently, the graph of the mapping $g(x)$ (and, by Lemma~\ref{lem:g}, the structure graph of the ring $\mathbb{Z}_p$) contains a cycle of length $>1$.
\end{proof}

The question of the existence of cycles of length $>1$ remains open for the case of primes of the form $p = 7 + 12k$.
The calculations carried out by the authors for $p \le 2077$ show that in this case too, cycles of length $>1$ always exist
(see Section \ref{sec:verif}).
Apparently, the following hypothesis is true.

\begin{hyp}\label{hyp:Zp}
The structure graph of the residue ring modulo a prime $\mathbb{Z}_p$ for $p \neq 2, 3, 7$ has $> p$ connected components.
\end{hyp}

\section{Cycles of Structure Graphs of Galois Fields of Extension Degree $>1$}

For Galois fields, an analogous hypothesis can be formulated.

\begin{hyp}\label{hyp:GF}
The structure graph of a Galois field of extension degree $>1$ different from $GF(4)$ contains a cycle of length $>1$.
\end{hyp}

For fields whose characteristic is not $2, 3, 7$, the validity of Hypothesis~\ref{hyp:GF} follows from Hypothesis~\ref{hyp:Zp}.

\begin{stmt}\label{stmt:GFpk}
If Hypothesis~\ref{hyp:Zp} is true, then the structure graph of the ring $GF(p^k)$ for $p \neq 2, 3, 7$ contains a cycle of length $>1$.
\end{stmt}
\begin{proof}
Indeed, $\mathbb{Z}_p \subseteq GF(p^k)$. From Hypothesis~\ref{hyp:Zp} it follows that in $\mathbb{Z}_p$ there is a cycle of length $>1$; it will also remain a cycle in $GF(p^k)$.
\end{proof}

\begin{stmt}\label{stmt:GF37}
The structure graph of rings of the form $GF(3^k)$ and $GF(7^k)$ for $k > 1$ contains a cycle of length $>1$.
\end{stmt}
\begin{proof}
The existence of cycles of length $>1$ in the rings $GF(3^2)$ and $GF(7^2)$ is established by direct verification (see Section \ref{sec:verif}).

Further, for even $k$ we have $GF(3^2) \subseteq GF(3^k)$, and a cycle in $GF(3^2)$ will be a cycle in $GF(3^k)$.

For odd $k$, elements of $GF(3^k) \setminus \mathbb{Z}_3$ cannot be mapped to $0$ (or, more generally, to any element of $\mathbb{Z}_3$).

Indeed, if $g(x) = c$ for $c \in \mathbb{Z}_3$, then $x$ is a root of the quadratic equation $x^2 - x + c = 0$ and, therefore, must lie either in $\mathbb{Z}_3$ or in $GF(3^2)$. But $GF(3^2) \cap GF(3^k) = \mathbb{Z}_3$. Thus, $x \notin GF(3^k) \setminus \mathbb{Z}_3$.

The assertion for $GF(7^k)$ is proved analogously.
\end{proof}

In the case of Galois fields of characteristic~2, we do not have the tool of passing to the graph of the mapping $g(x) = x - x^2$, and the existence of cycles of length $>1$ must be verified by computation. As shown above, it suffices to check only prime extension degrees and the extension degree~4 (since $GF(4)$ contains no cycles).

The authors have constructed the structure graphs of Galois fields of characteristic~2 up to $GF(2^{13})$. As a result of the computations, it was established that cycles of length $>1$ exist in the structure graphs of all rings except $GF(4)$ (see Section \ref{sec:verif}).

\section{Proof of Stages (3)--(5)}

Stage (3) follows from the existence of nilpotent elements in the matrix ring. By Statement~\ref{stmt:cycle2}, the structure graph contains a cycle of length~2.

Stage (4) follows from Statement~\ref{stmt:directsum}. Suppose that in at least one summand of the direct sum the structure graph contains a cycle of length $>1$. Take an element from this cycle and supplement it with arbitrary pairs of the form $(a, 0)$ from the other summands. As a result, we obtain a cycle of the same length in the direct sum.

Stage (5) follows from the nilpotency of the Jacobson radical (Statement~\ref{stmt:nilpotent}).

\section{Computational Verification}\label{sec:verif}

The computations mentioned in this paper can be reproduced using the scripts from the Zenodo repository \cite{link_repo}.
The computations were performed using the igraph Python package (\cite{link_igraph}).
For completeness and independence of exposition, the repository re-verifies the previously known cases.
See the table below for the mapping between statements, scripts, and results.

\begin{table}[ht]
\centering
\small
\caption{Mapping between statements, scripts, and results}
\label{Table_1}
\begin{tabular}{|p{0.4cm}|p{3.8cm}|p{1.7cm}|p{3.5cm}|p{2.6cm}|}
\hline
\textbf{\#} & \textbf{Statement} & \textbf{Prev. Publ.} & \textbf{Script} & \textbf{Results} \\
\hline
1 & $\mathbb{Z}_2$, $\mathbb{Z}_3$, $\mathbb{Z}_7$ \newline has $=n$ conn.comp. & \cite{link_Lipkovski_2012} & verify-small-rings.py & small-rings.csv\\
\hline
 2 & $GF(4)$ \newline has 4 conn.comp. & \cite{link_AS_2013}  & verify-gf4.py & gf4.csv\\
 \hline
 3 & $\mathbb{Z}_p$ for prime $p=7+12k$, \newline $19\leq p \leq 2077$ \newline has $>p$ conn.comp & for $ p < 200$ in \cite{link_Lipkovski_2015} & verify-primes-7-12k.py & primes-7-12k.csv\\
 \hline
 4 & $GF(9)$, $GF(49)$ \newline has $>n$ conn.comp. & --- & verify-gf9-gf49.py & gf9-gf49.csv \\
\hline
5 & $GF(2^m)$, $3\leq m \leq 7$  has $>2^m$ conn.comp. & --- & verify-gf8-gf128.py & gf8-gf128.csv\\
\hline
6 & $GF(2^{11})$ has $>2^{11}$ conn.comp. & --- & verify-gf2-11.py & gf2-11.csv\\
\hline
7 & $GF(2^{13})$, has $>2^{13}$ conn.comp. & --- & verify-gf2-13.py  & gf2-13.csv\\
\hline
\end{tabular}

\vspace{0.1cm}
\begin{flushleft}
\textit{The authors computations, \cite{link_repo}}
\end{flushleft}
\end{table}

\section{Conclusion}

In this article, the notion of the structure graph of an associative ring with unity is studied. The question is posed of classifying rings whose structure graph contains exactly $n$ connected components, which is equivalent to the absence of cycles of length $>1$.

The connection between the presence of cycles of length $>1$ in the structure graph of a ring and in the graph of the mapping $g(x) = x - x^2$ is shown. A hypothesis is formulated that the only rings whose structure graph contains no cycles of length $>1$ are the rings $\mathbb{Z}_2$, $\mathbb{Z}_3$, $\mathbb{Z}_7$, $GF(4)$ and their direct sums. The hypothesis is proved for all rings except residue rings modulo a prime of the form $7 + 12k$, $k \ge 181$, and Galois fields of characteristic~2 with extension degree $\ge 17$.

Further research can develop in the following directions: a final proof of this hypothesis, the study of cycle lengths of structure graphs, including the study of rings whose structure graph has cycles of large length.

\section*{Acknowledgments}

We are grateful to Professor Tuganbaev for his interest in our research and valuable discussions.

\vskip 0.4 true cm

\bibliographystyle{amsplain}

\end{document}